\documentclass{article}

\usepackage{latexsym}
\usepackage{amssymb}
\usepackage{amsthm}
\usepackage{amssymb,amsmath}
\usepackage{MnSymbol}
\usepackage{mathtools}
\usepackage{stmaryrd}
\usepackage{skak}
\usepackage{pigpen}
\usepackage[all]{xy}	
\usepackage{mathrsfs}
\usepackage{cancel}
\usepackage{harpoon}
\usepackage{yfonts}

\newcommand{\bpr}{\noindent{\em Proof\/. }}
\newcommand{\epr}{\hspace*{\fill}$\meddiamond$\medskip}

\newcommand{\bprLonebound}{\noindent{\em Proof of Lemma~\ref{L1bound}\/. }}

\newcommand{\bprthtwobound}{\noindent{\em Proof of Theorem~\ref{T2bound}\/. }}

\newtheorem{thm}{Theorem}
\newtheorem{lem}{Lemma}
\newtheorem{cor}{Corollary}
\newtheorem{prop}{Proposition}

\newtheorem{conj}{Conjecture}

\theoremstyle{definition}

\title{Determination of the prime bound of a graph}

\author{Abderrahim Boussa\"{\i}ri\thanks{Facult\'e des Sciences A\"{\i}n Chock, 
D\'epartement de Math\'ematiques et Informatique, Km 8 route d'El Jadida, 
BP 5366 Maarif, Casablanca, Maroc; {\tt aboussairi@hotmail.com}.}\and 
Pierre Ille\thanks{Institut de Math\'ematiques de Luminy, 
CNRS -- UMR 6206, 163 avenue de Luminy, Case~907, 13288 Marseille Cedex 09, France; {\tt ille@iml.univ-mrs.fr}.}
\thanks{Centre de recherches math\'ematiques, Universit\'e de Montr\'eal, Case postale 6128, 
Succursale Centre-ville, Montr\'eal, Qu\'ebec, Canada H3C 3J7.}}

\begin{document}

\maketitle

\begin{abstract}
Given a graph $G$, a subset $M$ of $V(G)$ is a module of $G$ if for each 
$v\in V(G)\setminus M$, $v$ is adjacent to all the elements of $M$ or to none of them. 
For instance, $V(G)$, $\emptyset$ and $\{v\}$ ($v\in V(G)$) are modules of $G$ called trivial. 
Given a graph $G$, $\omega_M(G)$ (respectively $\alpha_M(G)$) denotes the largest integer 
$m$ such that there is a module $M$ of $G$ which is a clique (respectively a stable) set in $G$ with  
$|M|=m$. 
A graph $G$ is prime if $|V(G)|\geq 4$ and if all its modules are trivial. 
The prime bound of $G$ is the smallest integer $p(G)$ such that there is a prime graph $H$ with $V(H)\supseteq V(G)$, 
$H[V(G)]=G$ and $|V(H)\setminus V(G)|=p(G)$. 
We establish the following. 
For every graph $G$ such that $\max(\alpha_M(G),\omega_M(G))\geq 2$ and 
$\log_2(\max(\alpha_M(G),\omega_M(G)))$ is not an integer, 
$p(G)=\lceil\log_2(\max(\alpha_M(G),\omega_M(G)))\rceil$. 
Then, we prove that for every graph $G$ such that $\max(\alpha_M(G),\omega_M(G))=2^k$ where $k\geq 1$, 
$p(G)=k$ or $k+1$. 
Moreover $p(G)=k+1$ if and only if $G$ or its complement admits $2^k$ isolated vertices. 
Lastly, we show that $p(G)=1$ for every non prime graph $G$ such that $|V(G)|\geq 4$ and 
$\alpha_M(G)=\omega_M(G)=1$. 
\end{abstract}

\medskip

\noindent {\bf Mathematics Subject Classifications (2010):}
05C70, 05C69

\medskip

\noindent {\bf Key words:} Module; prime graph; prime extension; prime bound; modular clique number; modular stability number

\section{Introduction}

A {\em graph} $G=(V(G),E(G))$ is constituted by a {\em vertex set} $V(G)$ and an {\em edge set} $E(G)\subseteq\binom{V(G)}{2}$. 
Given a set $S$, $K_S=(S,\binom{S}{2})$ is the {\em complete} graph on $S$ whereas $(S,\emptyset)$ is the {\em empty} graph. 
Let $G$ be a graph. 
With each $W\subseteq V(G)$ associate the {\em subgraph} $G[W]=(W,\binom{W}{2}\cap E(G))$ of $G$ induced by $W$. 
Given $W\subseteq V(G)$, $G[V(G)\setminus W]$ is also denoted by $G-W$ and by $G-w$ if $W=\{w\}$. 
A graph $H$ is an {\em extension} of $G$ if $V(H)\supseteq V(G)$ and $H[V(G)]=G$. 
Given $p\geq 0$, a $p$-extension of $G$ is an extension $H$ of $G$ such that $|V(H)\setminus V(G)|=p$. 
The {\em complement} of $G$ is the graph $\overline{G}=(V(G),\binom{V(G)}{2}\setminus E(G))$. 
A subset $W$ of $V(G)$ is a {\em clique} (respectively a {\em stable set}) in $G$ if $G[W]$ is complete (respectively empty). 
The largest cardinality of a clique (respectively a stable set) in $G$ is the {\em clique number} (respectively the {\em stability number}) of $G$, denoted by 
$\omega(G)$ (respectively $\alpha(G)$). 
Given $v\in V(G)$, the {\em neighbourhood} $N_G(v)$ of $v$ in $G$ is the family $\{w\in V(G):\{v,w\}\in E(G)\}$.
We consider $N_G$ as the function from $V(G)$ to $2^{V(G)}$ defined by $v\mapsto N_G(v)$ for each $v\in V(G)$. 
A vertex $v$ of $G$ is {\em isolated} if $N_G(v)=\emptyset$. 
The number of isolated vertices of $G$ is denoted by $\iota(G)$. 

We use the following notation. 
Let $G$ be a graph. 
For $v\neq w\in V(G)$, 
\begin{equation*}
(v,w)_G=
\begin{cases}
0&\text{if $\{v,w\}\not\in E(G)$},\\
1&\text{if $\{v,w\}\in E(G)$}. 
\end{cases}
\end{equation*}
Given $W\subsetneq V(G)$, $v\in V(G)\setminus W$ and $i\in\{0,1\}$, $(v,W)_G=i$ means $(v,w)_G=i$ for every $w\in W$. 
Given $W,W'\subsetneq V(G)$, with $W\cap W'=\emptyset$, and $i\in\{0,1\}$, $(W,W')_G=i$ means $(w,W')_G=i$ for every $w\in W$. 
Given $W\subsetneq V(G)$ and $v\in V(G)\setminus W$, $v\longleftrightarrow_GW$ means that there is $i\in\{0,1\}$ such that $(v,W)_G=i$. 
The negation is denoted by $v\not\longleftrightarrow_GW$.

Given a graph $G$, a subset $M$ of $V(G)$ is a {\em module} of $G$ if for each 
$v\in V(G)\setminus M$, we have $v\longleftrightarrow_GM$. 
For instance, $V(G)$, $\emptyset$ and $\{v\}$ ($v\in V(G)$) are modules of $G$ called {\em trivial}. 
Clearly, if $|V(G)|\leq 2$, then all the modules of $G$ are trivial. 
On the other hand, if $|V(G)|=3$, then $G$ admits a nontrivial module. 
A graph $G$ is then said to be {\em prime} if $|V(G)|\geq 4$ and if all its modules are trivial. 
For instance, given $n\geq 4$,  the {\em path} $(\{1,\ldots,n\},\{\{p,q\}:|p-q|=1\})$ is prime. 
Given a graph $G$, $G$ and $\overline{G}$ share the same modules. 
Thus $G$ is prime if and only if $\overline{G}$ is. 

Given a set $S$ with $|S|\geq 2$, 
$K_S$ admits a prime $\lceil\log_2(|S|+1)\rceil$-extension 
(see Sumner~\cite[Theorem~2.45]{S71} or Lemma~\ref{L1clique} below). 
This is extended to any graph in \cite[Theorem 3.7]{B07} and \cite[Theorem 3.2]{BRV11} as follows. 

\begin{thm}\label{Tbrignall}
A graph $G$, with $|V(G)|\geq 2$, admits a 
prime $\lceil\log_2(|V(G)|+1)\rceil$-extension. 
\end{thm}

Following Theorem~\ref{Tbrignall}, we introduce the notion of prime bound. 
Let $G$ be a graph. 
The {\em prime bound} of $G$ is the smallest integer $p(G)$ such that $G$ admits a prime 
$p(G)$-extension. 
Observe that $p(G)=p(\overline{G})$ for every graph $G$. 
By Theorem~\ref{Tbrignall}, $p(G)\leq\lceil\log_2(|V(G)|+1)\rceil$. 
By considering the clique number and the stability number, Brignall \cite[Conjecture 3.8]{B07} conjectured the following. 

\begin{conj}\label{Cbrignall}
For a graph $G$ with $|V(G)|\geq 2$, 
$$p(G)\leq\lceil\log_2(\max(\alpha(G),\omega(G))+1)\rceil.$$
\end{conj}

We answer the conjecture positively by refining the notions of clique number and of stability number as follows. 
Given a graph $G$, the {\em modular clique number} of $G$ is the largest integer 
$\omega_M(G)$ such that there is a module $M$ of $G$ which is a clique in $G$ with  
$|M|=\omega_M(G)$. 
The {\em modular stability number} of $G$ is 
$\alpha_M(G)=\omega_M(\overline{G})$. 
The following lower bound is simply obtained. 

\begin{lem}\label{L1bound}
For every  graph $G$ such that $\max(\alpha_M(G),\omega_M(G))\geq 2$, 
$$p(G)\geq\lceil\log_2(\max(\alpha_M(G),\omega_M(G)))\rceil.$$
\end{lem}

Theorem 3.2 of \cite{BRV11} is proved by induction on the number of vertices. 
Using the main arguments of this proof, 
we improve Theorem~\ref{Tbrignall} as follows. 

\begin{thm}\label{T1bound}
For every graph $G$ such that $\max(\alpha_M(G),\omega_M(G))\geq 2$, 
$$p(G)\leq\lceil\log_2(\max(\alpha_M(G),\omega_M(G))+1)\rceil.$$
\end{thm}

The proof of Theorem~\ref{T1bound} derives from an induction as well. 
A direct construction of a suitable extension is provided in~\cite[Theorem~2]{BI11}. 
The following is an immediate consequence of Lemma~\ref{L1bound} and Theorem~\ref{T1bound}. 

\begin{cor}\label{C1bound}
For every graph $G$ such that $\max(\alpha_M(G),\omega_M(G))\geq 2$, 
$$\lceil\log_2(\max(\alpha_M(G),\omega_M(G)))\rceil\leq p(G)\leq\lceil\log_2(\max(\alpha_M(G),\omega_M(G))+1)\rceil.$$
\end{cor}

Let $G$ be graph such that $\max(\alpha_M(G),\omega_M(G))\geq 2$. 
On the one hand, it follows from Corollary~\ref{C1bound} that 
\begin{equation*}
\max(\alpha_M(G),\omega_M(G))\not\in\{2^k:k\geq 1\}\ \Rightarrow\ p(G)=\lceil\log_2(\max(\alpha_M(G),\omega_M(G)))\rceil.
\end{equation*}
On the other, if $\max(\alpha_M(G),\omega_M(G))=2^k$, where $k\geq 1$, then 
$p(G)=k$ or $k+1$. 
The next allows us to determine this. 

\begin{thm}\label{T2bound}
For every graph $G$ such that $\max(\alpha_M(G),\omega_M(G))=2^k$ where $k\geq 1$, 
\begin{equation*}
\text{$p(G)=k+1$ if and only if $\iota(G)=2^k$ or 
$\iota(\overline{G})=2^k$.} 
\end{equation*}
\end{thm}

Lastly, we show that $p(G)=1$ for every non prime graph $G$ such that $|V(G)|\geq 4$ and 
$\alpha_M(G)=\omega_M(G)=1$ (see Proposition~\ref{m(G)=1}). 

\section{Preliminaries}

Given a graph $G$, the family of the modules of $G$ is denoted by $\mathcal{M}(G)$. 
Furthermore set $$\mathcal{M}_{\geq 2}(G)=\{M\in\mathcal{M}(G):|M|\geq 2\}.$$
We begin with the well known properties of the modules of a graph (for example, see \cite[Theorem~3.2, Lemma~3.9]{EHR99}). 

\begin{prop}\label{P1module}
Let $G$ be a graph. 
\begin{enumerate}
\item Given $W\subseteq V(G)$, 
$\ \{M\cap W:M\in\mathcal{M}(G)\}\subseteq\mathcal{M}(G[W])$. 
\item Given a module $M\in\mathcal{M}(G)$, 
$\ \mathcal{M}(G[M])=\{N\in\mathcal{M}(G):N\subseteq M\}$. 
\item Given $M,N\in\mathcal{M}(G)$ with $M\cap N=\emptyset$, there is $i\in\{0,1\}$ 
such that $(M,N)_G=i$.
\end{enumerate}
\end{prop}

Given a graph $G$, a partition $P$ of $V(G)$ is a {\em modular partition} of $G$ if 
$P\subseteq\mathcal{M}(G)$. 
Let $P$ be such a partition. 
Given $M\neq N\in P$, there is $i\in\{0,1\}$ such that $(M,N)_G=i$ by Proposition~\ref{P1module}.3. 
This justifies the following definition. 
The {\em quotient} of $G$ by $P$ is the graph $G/P$ defined on $V(G/P)=P$ by 
$(M,N)_{G/P}=(M,N)_G$ for $M\neq N\in P$. 
We use the following properties of the quotient (for example, see 
\cite[Theorems~4.1--4.3, Lemma~4.1]{EHR99}). 

\begin{prop}\label{P1quotient}
Given a graph $G$, consider a modular partition $P$ of $G$.
\begin{enumerate}
\item Given $W\subseteq V(G)$, if $|W\cap X|=1$ for each $X\in P$, then $G[W]$ and 
$G/P$ are isomorphic. 
\item For every $M\in\mathcal{M}(G)$, $\{X\in P:M\cap X\neq\emptyset\}\in\mathcal{M}(G/P)$. 
\item For every $Q\in\mathcal{M}(G/P)$, $\bigcup Q\in\mathcal{M}(G)$. 
\end{enumerate}
\end{prop}

The following strengthening of the notion of module is introduced to present the modular decomposition theorem (see Theorem~\ref{Tgallai} below). 
Given a graph $G$, a module $M$ of $G$ is said to be {\em strong} provided that for every $N\in\mathcal{M}(G)$, we have: if 
$M\cap N\neq\emptyset$, then $M\subseteq N$ or $N\subseteq M$. 
The family of the strong modules of $G$ is denoted by $\mathcal{S}(G)$. 
Furthermore set $$\mathcal{S}_{\geq 2}(G)=\{M\in\mathcal{S}(G):|M|\geq 2\}.$$
We recall the following well known properties of the strong modules of a 
graph (for example, see 
\cite[Theorem~3.3]{EHR99}).

\begin{prop}\label{P1strong}
Let $G$ be a graph. 
For every $M\in\mathcal{M}(G)$, $\mathcal{S}(G[M])=\{N\in\mathcal{S}(G):N\subsetneq M\}\cup\{M\}$.
\end{prop}

With each graph $G$, we associate the family $\Pi(G)$ of the maximal proper and nonempty strong modules of $G$ under inclusion. 
For convenience set $$\Pi_1(G)=\{M\in\Pi(G):|M|=1\}\ \text{and}\ \Pi_{\geq 2}(G)=\{M\in\Pi(G):|M|\geq 2\}.$$
The modular decomposition theorem is stated as follows. 

\begin{thm}[Gallai \cite{G67,MP01}]\label{Tgallai}
For a graph $G$ with $|V(G)|\geq 2$, the family $\Pi(G)$ realizes a modular partition of $G$. 
Moreover, the corresponding quotient $G/\Pi(G)$ is complete, empty or prime. 
\end{thm}

Let $G$ be a graph with $|V(G)|\geq 2$. 
As a direct consequence of the definition of a strong module, we obtain that 
the family $\mathcal{S}(G)\setminus\{\emptyset\}$ endowed with inclusion is a tree called the {\em modular decomposition tree} \cite{CM05} of $G$. 
Given $M\in\mathcal{S}_{\geq 2}(G)$, it follows from 
Proposition~\ref{P1strong} that $\Pi(G[M])\subseteq\mathcal{S}(G)$. 
Furthermore, given $W\subseteq V(G)$, 
the family $\{M\in\mathcal{S}(G):M\supseteq W\}$ endowed with inclusion is a total order. 
Its smallest element is denoted by $\widehat{W}$.

Let $G$ be a graph with $|V(G)|\geq 2$. 
Using Theorem~\ref{Tgallai}, we label $\mathcal{S}_{\geq 2}(G)$ by the function $\lambda_G$ defined as follows. 
For each $M\in\mathcal{S}_{\geq 2}(G)$, 
\begin{equation*}
\lambda_G(M)=
\begin{cases}
\filledmedsquare\ \ \text{if $G[M]/\Pi(G[M])$ is complete,}\\
\medsquare\ \ \text{if $G[M]/\Pi(G[M])$ is empty,}\\
\sqcup\ \ \text{if $G[M]/\Pi(G[M])$ is prime.}
\end{cases}
\end{equation*}

\section{Some prime extensions}

\begin{lem}\label{L1stable}
Let $S$ and $S'$ be disjoint sets such that $|S|\geq 2$ and $|S'|=\lceil\log_2(|S|+1)\rceil$. 
There exists a prime graph $G$ defined on $V(G)=S\cup S'$ such that $S$ and $S'$ are stable sets in $G$. 
\end{lem}

\bpr
If $|S|=2$, then $|S'|=2$ and we can choose a path on 4 vertices for $G$. 
Assume that $|S|\geq 3$. 
As $|S'|=\lceil\log_2(|S|+1)\rceil$, $2^{|S'|-1}\leq |S|$ and hence 
$|S'|\leq|S|$. 
Thus there exists a bijection $\psi_{S'}$ from $S'$ onto $S''\subseteq S$. 
Consider the injection $f_{S''}:S''\longrightarrow 2^{S'}\setminus\{\emptyset\}$ defined by 
$s''\mapsto S'\setminus\{(\psi_{S'})^{-1}(s'')\}$. 
Since $|S'|=\lceil\log_2(|S|+1)\rceil$, $|S|<2^{|S'|}$ and there exists an injection $f_S$ from $S$ into $2^{S'}\setminus\{\emptyset\}$ such that 
$(f_S)_{\restriction S''}=f_{S''}$. 
Lastly, consider the graph $G$ defined on $V(G)=S\cup S'$ such that $S$ and $S'$ are stable sets  in $G$ and $(N_G)_{\restriction S}=f_S$. 
We prove that $G$ is prime. 
If $|S|=3$, then $|S'|=2$ and $G$ is a path on 5 vertices which is prime. 
Assume that $|S|\geq 4$ and hence $|S'|\geq 3$. 
Let $M\in\mathcal{M}_{\geq 2}(G)$. 

First, if $M\subseteq S$, then we would have $f_S(u)=f_S(v)$ for any $u\neq v\in M$. 
Thus $M\cap S'\neq\emptyset$. 

Second, suppose that $M\subseteq S'$. 
Recall that for each $s\in S$, either $M\cap N_G(s)=\emptyset$ or $M\subseteq N_G(s)$. 
Given $u\in M$, 
consider the function $f:S\longrightarrow 2^{((S'\setminus M)\cup\{u\})}\setminus\{\emptyset\}$ defined by 
\begin{equation*}
f(s)=
\begin{cases}
\text{$N_G(s)$ if $M\cap N_G(s)=\emptyset$,}\\
\text{$(N_G(s)\setminus M)\cup\{u\}$ if $M\subseteq N_G(s)$,}
\end{cases}
\end{equation*}
for every $s\in S$. 
Since $(N_G)_{\restriction S}$ is injective, $f$ is also and we would obtain that 
$\left|S\right|<2^{\left|S'\right|-1}$. 
It follows that $M\cap S\neq\emptyset$. 

Third, suppose that $S'\setminus M\neq\emptyset$. 
We have $(S\cap M,S'\setminus M)_G=(S'\cap M,S'\setminus M)_G=0$. 
Given $s'\in S'\cap M$, $N_G(\psi_{S'}(s'))=S'\setminus\{s'\}$. 
In particular $S'\setminus M\subseteq N_G(\psi_{S'}(s'))$ and hence $\psi_{S'}(s')\in S\setminus M$. 
Furthermore $(\psi_{S'}(s'),S'\cap M)_G=(\psi_{S'}(s'),S\cap M)_G=0$. 
Therefore $S'\cap M=\{s'\}$. 
Similarly, we prove that $|S'\setminus M|=1$ which would imply that $|S'|=2$. 
It follows that $S'\subseteq M$. 

Lastly, suppose that $S\setminus M\neq\emptyset$. 
For each $s\in S\setminus M\neq\emptyset$, we would have $(s,S')_G=(s,S\cap M)_G=0$ and hence 
$N_G(s)=\emptyset$. 
It follows that $S\subseteq M$ and $M=S\cup S'$. 
\epr

\begin{lem}\label{L1clique}
Let $C$ and $S'$ be disjoint sets such that $|C|\geq 2$ and $|S'|=\lceil\log_2(|C|+1)\rceil$. 
There exists a prime graph $G$ defined on $V(G)=C\cup S'$ such that $C$ is a clique and $S'$ is a stable set in $G$. 
\end{lem}

\bpr
There exists a bijection $\psi_{S'}$ from $S'$ onto $S''\subseteq C$. 
Consider the injection $f_{S''}:S''\longrightarrow 2^{S'}\setminus\{S'\}$ defined by 
$s''\mapsto\{(\psi_{S'})^{-1}(s'')\}$. 
Let $f_S$ be any injection from $S$ into $2^{S'}\setminus\{S'\}$ such that 
$(f_S)_{\restriction S''}=f_{S''}$. 
Lastly, consider the graph $G$ defined on $V(G)=C\cup S'$ such that $C$ is a clique in $G$, $S'$ is a stable set in $G$ and 
$N_G(c)\cap S'=f_S(c)$ for each $c\in C$. 
We prove that $G$ is prime. 
Let $M\in\mathcal{M}_{\geq 2}(G)$. 
As in the proof of Lemma~\ref{L1stable}, we have $M\cap C\neq\emptyset$ and $M\cap S'\neq\emptyset$. 

Now, suppose that $S'\setminus M\neq\emptyset$. 
We have $(C\cap M,S'\setminus M)_G=(S'\cap M,S'\setminus M)_G=0$. 
Given $t'\in S'\setminus M$, $N_G(\psi_{S'}(t'))\cap S'=\{t'\}$. 
Thus $\psi_{S'}(t')\in C\setminus M$. 
But $(\psi_{S'}(t'),S'\cap M)_G=(\psi_{S'}(t'),C\cap M)_G=1$ which contradicts $N_G(\psi_{S'}(t'))\cap S'=\{t'\}$. 
It follows that $S'\subseteq M$. 

Lastly, suppose that $C\setminus M\neq\emptyset$. 
For each $c\in C\setminus M\neq\emptyset$, we would have $(c,S')_G=(c,C\cap M)_G=1$ and hence 
$N_G(c)\cap S'=S'$. 
It follows that $S\subseteq M$ and $M=S\cup S'$. 
\epr

The question of prime extensions of a prime graph is not detailed enough in \cite{BRV11}. 
For instance, the number of prime 1-extensions of a prime graph given in \cite{BRV11} is not correct. 
Moreover, Corollary~\ref{C1extension} below is used without a precise proof. 

\begin{lem}\label{L1extension}
Let $G$ be a prime graph $G$. Given $a\not\in V(G)$, there exist 
$$2^{|V(G)|}-2|V(G)|-2$$ distinct prime extensions of $G$ to $V(G)\cup\{a\}$. 
\end{lem}

\bpr
Consider any graph $H$ defined on $V(H)=V(G)\cup\{a\}$ such that $H[V(G)]=G$. 
We prove that $H$ is not prime if and only if 
$$N_H(a)\in\{\emptyset,V(G)\}\cup\{N_G(v):v\in V(G)\}\cup\{N_G(v)\cup\{v\}:v\in V(G)\}.$$
To begin, assume that 
$N_H(a)\in \{\emptyset,V(G)\}\cup\{N_G(v):v\in V(G)\}\cup\{N_G(v)\cup\{v\}:v\in V(G)\}$. 
If $N_H(a)=\emptyset$ or $V(G)$, then $V(G)$ is a nontrivial module of $H$. 
If there is $v\in V(G)$ such that $N_H(a)\setminus\{v\}=N_G(v)$, then $\{a,v\}$ is a nontrivial module of $H$. 

Conversely, assume that $H$ admits a nontrivial module $M$. 
By Proposition~\ref{P1module}.1, $M\setminus\{a\}\in\mathcal{M}(G)$. 
As $G$ is prime and as $M\setminus\{a\}\neq\emptyset$ and $M\subsetneq V(H)$, either $|M\setminus\{a\}|=1$ or $M=V(G)$. 
In the second instance, $N_H(a)=\emptyset$ or $V(G)$. 
In the first, there is $v\in V(G)$ such that $M=\{a,v\}$. 
Thus $N_H(a)=N_G(v)$ or $N_G(v)\cup\{v\}$. 

To conclude, observe that $$|\{\emptyset,V(G)\}\cup\{N_G(v):v\in V(G)\}\cup\{N_G(v)\cup\{v\}:v\in V(G)\}|=2+2|V(G)|$$ because $G$ is prime. 
\epr

\begin{cor}\label{C1extension}
Let $G$ be a prime graph $G$. For any $a\neq b\not\in V(G)$, there exists a prime extension $H$ of $G$ to $V(G)\cup\{a,b\}$ such that $(a,b)_H=0$. 
\end{cor}

\bpr
Since $|V(G)|\geq 4$, $2^{|V(G)|}-2|V(G)|-2\geq 2$. 
Consequently there is an extension $H$ of $G$ to $V(G)\cup\{a,b\}$ such that 
$(a,b)_H=0$, $N_H(a)\neq N_H(b)$ and 
$$N_H(a), N_H(b)\not\in \{\emptyset,V(G)\}\cup\{N_G(v):v\in V(G)\}\cup\{N_G(v)\cup\{v\}:v\in V(G)\}.$$
By the proof of Lemma~\ref{L1extension}, $H-a$ and $H-b$ are prime. 
We show that $H$ is prime also. 
Let $M\in\mathcal{M}_{\geq 2}(H)$. 
By Proposition~\ref{P1module}.1, $M\setminus\{a\}\in\mathcal{M}(H-a)$. 
As $H-a$ is prime and $M\setminus\{a\}\neq\emptyset$, either $|M\setminus\{a\}|=1$ or 
$M\setminus\{a\}=V(H)\setminus\{a\}$. 
In the first, there is $v\in V(G)\cup\{b\}$ such that $M=\{a,v\}$. 
If $v=b$, then $N_H(a)=N_H(b)$. 
If $v\in V(G)$, then $\{a,v\}$ would be a nontrivial module of $H-b$. 
Consequently $M\setminus\{a\}=V(H)\setminus\{a\}$. 
Since $H-b$ is prime, $a\not\longleftrightarrow_H V(G)$ and hence $a\in M$. 
Thus $M=V(H)$. 
\epr

\section{Proof of Theorem~\ref{T1bound}}

Let $G$ be a graph with $|V(G)|\geq 2$. 
By \cite[Theorem 3.2]{BRV11}, there exists a prime extension $H$ of $G$ such that 
\begin{equation*}
\begin{cases}
2\leq|V(H)\setminus V(G)|\leq\lceil\log_2(|V(G)|+1)\rceil\\
\text{and}\\
\text{$V(H)\setminus V(G)$ is a stable set in $H$.}
\end{cases}
\end{equation*}
We can consider the smallest integer $q(G)$ such that $q(G)\geq 2$ and $G$ admits a prime 
$q(G)$-extension $H$ such that $V(H)\setminus V(G)$ is a stable set in $H$. 

The results below, from Proposition~\ref{P1th1} to Corollary~\ref{C2th1}, are suggested by the proof of \cite[Theorem 3.2]{BRV11}. 

We introduce a basic construction. 
Consider a graph $G$ and a modular partition $P$ of $G$ such that $P\subseteq\mathcal{S}(G)$ and 
$P\cap\mathcal{S}_{\geq 2}(G)\neq\emptyset$. 
Let $X\in P\cap\mathcal{S}_{\geq 2}(G)$ such that 
$$q(G[X])=\max(\{q(G[Y]):Y\in P\cap\mathcal{S}_{\geq 2}(G)\}).$$ 
Consider a set $S$ such that $S\cap V(G)=\emptyset$ and $|S|=q(G[X])$. 
There exists a prime $q(G[X])$-extension $H_X$ of $G[X]$ to $X\cup S$ such that $S$ is a stable set in $H_X$. 
Since $X$ is not a module of $H_X$, there is $s_X\in S$ such that $s_X\not\longleftrightarrow_G X$. 
Furthermore, if there is $v\in S$ such that $(v,X)_{H_X}=0$, then 
$V(H_X)\setminus\{v\}$ would be a nontrivial module of $H_X$. 
Thus $\{v\in S:v\longleftrightarrow_{H_X} X\}=\{v\in S:(v,X)_{H_X}=1\}$. 
As $S$ is a stable set in $H_X$, $\{v\in S:(v,X)_{H_X}=1\}$ is a module of $G$. 
It follows that 
\begin{equation*}
\begin{cases}
\{v\in S:v\longleftrightarrow_{H_X} X\}=\{v\in S:(v,X)_{H_X}=1\}\\
\text{and}\\
|\{v\in S:v\longleftrightarrow_{H_X} X\}|\leq 1\\
\text{and}\\
s_X\in S\setminus\{v\in S:v\longleftrightarrow_{H_X} X\}.
\end{cases}
\end{equation*}
Now, for each $Y\in (P\cap\mathcal{S}_{\geq 2}(G))\setminus\{X\}$, there is a prime $q(G[Y])$-extension $H_Y$ of 
$G[Y]$ to $Y\cup S_Y$ such that $\{v\in S:v\longleftrightarrow_{H_X} X\}\subseteq S_Y\subseteq S$ and $S_Y$ is a stable set in $H_Y$. 
Consider the extension $H$ of $G$ to $V(G)\cup S$ satisfying
\begin{itemize}
\item for each $Y\in P\cap\mathcal{S}_{\geq 2}(G)$, $H[Y\cup S_Y]=H_Y$;
\item for each $v\in V(G)$ such that $\{v\}\in P$, $(v,S\setminus\{s_X\})_H=0$ and 
$(v,s_X)_H=1$. 
\end{itemize}

\begin{prop}\label{P1th1}
Given a graph $G$, consider a modular partition $P$ of $G$ such that 
$P\subseteq\mathcal{S}(G)$ and $P\cap\mathcal{S}_{\geq 2}(G)\neq\emptyset$. 
If the corresponding extension $H$ is not prime, then all the nontrivial modules of $G$ are included in $\{v\in V(G):\{v\}\in P\}$.
\end{prop}

\bpr
Let $M$ be a nontrivial module of $H$. 
By Proposition~\ref{P1module}.1, $M\cap(X\cup S)\in\mathcal{M}(H[X\cup S])$. 
Since $H[X\cup S]$ is prime, we have $M\supseteq X\cup S$, $|M\cap(X\cup S)|=1$ or $M\cap(X\cup S)=\emptyset$. 

For a first contradiction, suppose that $M\supseteq X\cup S$. 
Let $v\in V$ such that $\{v\}\in P$. 
As $v\not\longleftrightarrow_H S$, $v\in M$. 
Thus $\{v\in V(G):\{v\}\in P\}\subseteq M$. 
Let $Y\in P\cap\mathcal{S}_{\geq 2}(G)$. 
By Proposition~\ref{P1module}.1, $M\cap(Y\cup S_Y)\in\mathcal{M}(H[Y\cup S_Y])$. 
Since $H[Y\cup S_Y]$ is prime and since $S_Y\subseteq M\cap(Y\cup S_Y)$, $Y\subseteq M$. 
Therefore $\bigcup(P\cap\mathcal{S}_{\geq 2}(G))\subseteq M$ and we would have $M=V(H)$. 

For a second contradiction, suppose that $|M\cap(X\cup S)|=1$. 
Consider $v\in S\cup X$ such that $M\cap(X\cup S)=\{v\}$. 
Suppose that $v\in X$. 
We have $M\subseteq V(G)$ and $M\in\mathcal{M}(G)$ by Proposition~\ref{P1module}.1. 
As $X\in\mathcal{S}(G)$ and $v\in X\cap M$, $X\subseteq M$ or $M\subseteq X$. 
In both cases, we would have $|M\cap(X\cup S)|\geq 2$. 
Suppose that $v\in S$. 
There is $Y\in P\setminus\{X\}$ such that $Y\cap M\neq\emptyset$. 
Let $y\in Y\cap M$. 
Since $y\longleftrightarrow_G X$, $v\longleftrightarrow_{H_X} X$ and hence $v\neq s_X$. 
If $Y\in P\cap\mathcal{S}_{\geq 2}(G)$, then $v\in S_Y$ and $M\cap (Y\cup S_Y)$ would be a nontrivial module of $H[Y\cup S_Y]$. 
If $Y=\{y\}$, then $(y,s_X)_H=1$. 
Thus $(v,s_X)_H=1$ and $S$ would not be a stable set in $H$. 

It follows that $M\cap(X\cup S)=\emptyset$. 
By Proposition~\ref{P1module}.1, $M\in\mathcal{M}(G)$. 
Let $Y\in (P\cap\mathcal{S}_{\geq 2}(G))\setminus\{X\}$. 
Suppose for a contradiction that $Y\cap M\neq\emptyset$. 
As $Y\in\mathcal{S}(G)$, $Y\subseteq M$ or $M\subseteq Y$. 
In both cases, $M\cap (Y\cup S_Y)$ would be a nontrivial module of $H[Y\cup S_Y]$. 
It follows that $Y\cap M=\emptyset$. 
Therefore $M\subseteq\{v\in V(G):\{v\}\in P\}$. 
\epr

\begin{cor}\label{C1th1}
Given a graph $G$ such that $G/\Pi(G)$ is prime, we have 
\begin{equation*}
q(G)\leq\ 
\begin{cases}
2\ \text{if}\ \ \Pi_{\geq 2}(G)=\emptyset\\
\max(\{q(G[X]):X\in\Pi_{\geq 2}(G)\})\ \text{if}\ \ \Pi_{\geq 2}(G)\neq\emptyset.
\end{cases}
\end{equation*}
\end{cor}

\bpr
If $G$ is prime, then $q(G)\leq 2$ by Corollary~\ref{C1extension}, and hence $q(G)=2$. 
Assume that $G$ is not prime, that is, $\Pi_{\geq 2}(G)\neq\emptyset$. 
Let $H$ be the extension of $G$ associated with $\Pi(G)$. 
Suppose that $H$ admits a nontrivial module $M$. 
By Proposition~\ref{P1th1}, $\{\{u\}:u\in M\}\subseteq\Pi_1(G)$. 
Thus $M\in\mathcal{M}(G)$ by Proposition~\ref{P1module}.1. 
By Proposition~\ref{P1quotient}.2, $\{\{u\}:u\in M\}$ would be a nontrivial module of $G/\Pi(G)$. 
\epr

\begin{prop}\label{P2th1}
Given a graph $G$ such that $G/\Pi(G)$ is complete or empty, we have 
\begin{equation*}
\begin{cases}
q(G)\leq\max(2,\lceil\log_2(|\Pi_1(G)|+1)\rceil)\\
\text{or}\\
q(G)\leq\max(\{q(G[X]):X\in\Pi_{\geq 2}(G)\}).
\end{cases}
\end{equation*}
\end{prop}

\bpr
Assume that $G/\Pi(G)$ is empty. 
If $\Pi(G)=\Pi_1(G)$, then $G$ is empty by Proposition~\ref{P1quotient}.1, and it suffices to apply Lemma~\ref{L1stable}. 
Assume that $\Pi_{\geq 2}(G)\neq\emptyset$ and set $$W_2=\bigcup\Pi_{\geq 2}(G).$$ 
Let $H$ be the extension of $G$ associated with $\Pi(G)$. 
Recall that $V(H)=V(G)\cup S$, $V(G)\cap S=\emptyset$ and $|S|=q(G[X])$ where $X\in\Pi_{\geq 2}(G)$ such that 
$q(G[X])=\max(\{q(G[Y]):Y\in\Pi_{\geq 2}(G)\})$. 
Moreover $H[X\cup S]$ is prime. 

If $|\Pi_1(G)|\leq 1$, then $H$ is prime by Proposition~\ref{P1th1} so that $q(G)\leq\max(\{q(G[Y]):Y\in\Pi_{\geq 2}(G)\})$. 
Assume that $|\Pi_1(G)|\geq 2$ and set $$W_1=V(G)\setminus W_2.$$ 
By Lemma~\ref{L1stable}, there exists a prime extension $H_1$ of $G[W_1]$ to $W_1\cup S_1$ such that 
$|S_1|=\lceil\log_2(|W_1|+1)\rceil$ and $S_1$ is stable in $H_1$. 
As $G/\Pi(G)$ is empty, $\Pi_{\geq 2}(G)\in\mathcal{M}(G/\Pi(G))$. 
By Proposition~\ref{P1quotient}.3, $W_2\in\mathcal{M}(G)$. 
Thus $\Pi_{\geq 2}(G)\subseteq\mathcal{S}(G[W_2])$ by Proposition~\ref{P1strong}. 
It follows from Proposition~\ref{P1th1} that $H[W_2\cup S]$ is prime. 
We construct suitable extensions of $G$ according to whether $|S_1|\leq|S|$ or not. 

To begin, assume that $|S_1|\leq|S|$. 
We can assume that $$\{v\in S:v\longleftrightarrow_{H[X\cup S]}X\}\subseteq S_1\subseteq S$$ and we consider an extension $H'$ of $H_1$ and $H[W_2\cup S]$ to $V(G)\cup S$. 
We show that $H'$ is prime. 
Let $M\in\mathcal{M}_{\geq 2}(H')$. 
By Proposition~\ref{P1module}.1, $M\cap(W_2\cup S)\in\mathcal{M}(H[W_2\cup S])$. 
Since $H[W_2\cup S]$ is prime, $M\cap(W_2\cup S)=\emptyset$, $|M\cap(W_2\cup S)|=1$ or $M\supseteq(W_2\cup S)$. 
\begin{itemize}
\item Suppose for a contradiction that $M\cap(W_2\cup S)=\emptyset$. 
By Proposition~\ref{P1module}.1, $M$ would be a nontrivial module of $H_1$. 
\item Suppose for a contradiction that $|M\cap(W_2\cup S)|=1$ and consider $w\in W_2\cup S$ such that $M\cap(W_2\cup S)=\{w\}$. 
First, suppose that $w\in W_2$ and consider $Y\in\Pi_{\geq 2}(G)$ such that $w\in Y$. 
By Proposition~\ref{P1module}.1, $M\in\mathcal{M}(G)$. 
As $Y\in\mathcal{S}(G)$ and $w\in X\cap M$, $X\subseteq M$ or $M\subseteq X$. 
In both cases, we would have $|M\cap(W_2\cup S)|\geq 2$. 
Second, suppose that $w\in S$ and consider $v\in W_1\cap M$. 
Since $v\longleftrightarrow_G X$, $w\longleftrightarrow_{H[W_2\cup S]}X$ and hence $w\in S_1$. 
It follows from Proposition~\ref{P1module}.1 that $M$ would be a nontrivial module of $H_1$. 
\end{itemize}
Consequently $M\supseteq(W_2\cup S)$. 
By Proposition~\ref{P1module}.1, $M\cap(W_1\cup S_1)\in\mathcal{M}(H_1)$. 
As $H_1$ is prime and $M\cap(W_1\cup S_1)\supseteq S_1$, $M\cap(W_1\cup S_1)=(W_1\cup S_1)$ so that $M=V(H')$. 

Now, assume that $|S_1|>|S|$. 
We can assume that $S\subsetneq S_1$ and we consider the unique extension $H''$ of $H_1$ and $H[W_2\cup S]$ to $V(G)\cup S_1$ such that 
\begin{equation}\label{E1P2th1}
(W_2,S_1\setminus S)_{H''}=0.
\end{equation}
We show that $H''$ is prime. 
Let $M\in\mathcal{M}_{\geq 2}(H'')$. 
We obtain $M\cap(W_1\cup S_1)=\emptyset$, $|M\cap(W_1\cup S_1)|=1$ or $M\supseteq(W_1\cup S_1)$. 
If $M\cap(W_1\cup S_1)=\emptyset$, then $M$ would be a nontrivial module of $H[W_2\cup S]$. 

Suppose for a contradiction that $|M\cap(W_1\cup S)_1|=1$ and consider $w\in W_1\cup S_1$ such that $M\cap(W_1\cup S_1)=\{w\}$. 
There is $v\in W_2\cap M$. 
Let $Y\in\Pi_{\geq 2}(G)$ such that $v\in Y$. 
\begin{itemize}
\item Suppose that $w\in W_1$. 
By Proposition~\ref{P1module}.1, $M\in\mathcal{M}(G)$. 
Since $Y\in\mathcal{S}(G)$ and since $Y\cap M\neq\emptyset$ and $w\in M\setminus Y$, $Y\subseteq M$. 
It follows from Proposition~\ref{P1module}.1 that $M\cap(W_2\cup S)$ would be a nontrivial module of $H[W_2\cup S]$. 
\item Suppose that $w\in S_1$. 
By Proposition~\ref{P1module}.1, $M\cap(W_2\cup S)\in\mathcal{M}(H[W_2\cup S])$. 
As $H[W_2\cup S]$ is prime and as $v\in M\cap W_2$ and $M\cap S\subseteq\{w\}$, $M\cap(W_2\cup S)=\{v\}$ and hence 
$w\in S_1\setminus S$. 
For every $u\in W_2\setminus\{v\}$, we have $(u,v)_G=(u,w)_{H''}=0$ by \eqref{E1P2th1}. 
Since $(v,W_1)_G=0$, we would have $N_G(v)=\emptyset$ and hence $\{v\}\in\Pi_1(G)$. 
\end{itemize}
It follows that $M\supseteq(W_1\cup S_1)$. 
By Proposition~\ref{P1module}.1, $M\cap(W_2\cup S)\in\mathcal{M}(H[W_2\cup S])$. 
As $H[W_2\cup S]$ is prime and $M\cap(W_2\cup S)\supseteq S$, $M\cap(W_2\cup S)=(W_2\cup S)$ so that $M=V(H'')$. 

Finally, observe that when $G/\Pi(G)$ is complete, we can proceed as previously by replacing \eqref{E1P2th1} by 
$(W_2,S_1\setminus S)_{H''}=1$. 
\epr

The next result follows from Corollary~\ref{C1th1} and Proposition~\ref{P2th1} by climbing the modular decomposition tree from bottom to top. 

\begin{cor}\label{C2th1}
Given a graph $G$, if there is $X\in\mathcal{S}_{\geq 2}(G)$ such that $\lambda_G(X)\in\{\medsquare,\filledmedsquare\}$ and 
$|\Pi_1(G[X])|\geq 2$, then 
\begin{equation*}
q(G)\leq\max(\{\lceil\log_2(|\Pi_1(G[Y])|+1)\rceil:Y\in\mathcal{S}_{\geq 2}(G),\lambda_G(Y)\in\{\medsquare,\filledmedsquare\}\}).
\end{equation*}
\end{cor}

Given Corollary~\ref{C2th1}, Theorem~\ref{T1bound} follows from the next transcription in terms of the modular decomposition tree. 
Let $G$ be a graph. 
Denote by $\mathbb{M}(G)$ the family of the maximal elements of $\mathcal{M}_{\geq 2}(G)$ under inclusion which are cliques or stable sets in $G$. 

\begin{prop}\label{P3th1}
Given a graph $G$ such that $\max(\alpha_M(G),\omega_M(G))\geq 2$, 
\begin{equation*}
M\in\mathbb{M}(G)\ \iff\ 
\begin{cases}
M\in\mathcal{M}_{\geq 2}(G)\\
\text{and}\\
\lambda_G(\widehat{M})\in\{\medsquare,\filledmedsquare\}\\
\text{and}\\
M=\{v\in\widehat{M}:\{v\}\in\Pi(G[\widehat{M}])\}.
\end{cases}
\end{equation*}
\end{prop}

\bpr
To begin, consider $M\in\mathbb{M}(G)$ and assume that $M$ is a stable set in $G$. 
By Proposition~\ref{P1module}.1, $M\in\mathcal{M}(G[\widehat{M}])$. 
Set $$Q=\{X\in\Pi_1(G[\widehat{M}]):X\cap M\neq\emptyset\}.$$ 
By definition of $\widehat{M}$, $|Q|\geq 2$ and hence $M=\bigcup Q$ because $Q\subseteq\mathcal{S}(G[\widehat{M}])$. 
Furthermore, $Q\subseteq\mathcal{S}(G[M])$ by Proposition~\ref{P1strong}. 
As all the strong modules of an empty graph are trivial, we obtain $|X|=1$ for each $X\in Q$, that is, 
$$M\subseteq\{v\in\widehat{M}:\{v\}\in\Pi(G[\widehat{M}])\}.$$ 
By Proposition~\ref{P1quotient}.2, $Q\in\mathcal{M}(G[\widehat{M}]/\Pi(G[\widehat{M}]))$. 
For a contradiction, suppose that $\lambda_G(\widehat{M})=\sqcup$. 
Since $Q\in\mathcal{M}_{\geq 2}(G[\widehat{M}]/\Pi(G[\widehat{M}]))$, $Q=\Pi(G[\widehat{M}])$ and hence $M=\widehat{M}$. 
As $|X|=1$ for each $X\in Q$, $G[\widehat{M}]/\Pi(G[\widehat{M}])$ and $G[\widehat{M}]$ are isomorphic by Proposition~\ref{P1quotient}.1. 
It would follow that $G[M]$ is prime. 
Consequently $\lambda_G(\widehat{M})\in\{\medsquare,\filledmedsquare\}$. 
Given $v\neq w\in M$, we have $(\{v\},\{w\})_{G[\widehat{M}]/\Pi(G[\widehat{M}])}=(v,w)_G=0$. 
Thus $$\lambda_G(\widehat{M})=\medsquare.$$
Since $\lambda_G(\widehat{M})=\medsquare$, $\Pi_1(G[\widehat{M}])\in\mathcal{M}(G[\widehat{M}]/\Pi(G[\widehat{M}]))$. 
By Proposition~\ref{P1quotient}.3, $\bigcup\Pi_1(G[\widehat{M}])\in\mathcal{M}(G[\widehat{M}])$ and hence 
$\bigcup\Pi_1(G[\widehat{M}])\in\mathcal{M}(G)$ by Proposition~\ref{P1module}.2. 
Given $v\neq w\in\bigcup\Pi_1(G[\widehat{M}])$, we have $(v,w)_G=(\{v\},\{w\})_{G[\widehat{M}]/\Pi(G[\widehat{M}])}=0$. 
Therefore $\bigcup\Pi_1(G[\widehat{M}])$ is a stable set of $G$. 
As $M\subseteq\bigcup\Pi_1(G[\widehat{M}])$, $M=\bigcup\Pi_1(G[\widehat{M}])$ by maximality of $M$. 
It follows that $$M=\{v\in\widehat{M}:\{v\}\in\Pi(G[\widehat{M}])\}.$$ 

Conversely, consider $M\in\mathcal{M}_{\geq 2}(G)$ such that $\lambda_G(\widehat{M})=\medsquare$ and 
$M=\{v\in\widehat{M}:\{v\}\in\Pi(G[\widehat{M}])\}$. 
As $\lambda_G(\widehat{M})=\medsquare$, 
$\Pi_1(G[\widehat{M}])\in\mathcal{M}(G[\widehat{M}]/\Pi(G[\widehat{M}]))$. 
By Proposition~\ref{P1quotient}.3, 
$M=\bigcup\Pi_1(G[\widehat{M}])\in\mathcal{M}(G[\widehat{M}])$ and hence 
$M\in\mathcal{M}(G)$ by Proposition~\ref{P1module}.2. 
Since $(v,w)_G=(\{v\},\{w\}_{G[\widehat{M}]/\Pi(G[\widehat{M}])}=0$ for $v\neq w\in M$, 
$M$ is a stable set in $G$. 
There is $N\in\mathbb{M}(G)$ such that $N\supseteq M$. 
As $M$ is a stable set in $G$, $N$ is as well. 
By what precedes, $N=\{v\in\widehat{N}:\{v\}\in\Pi(G[\widehat{N}])\}$. 
We have $\widehat{M}\subseteq\widehat{N}$ because $M\subseteq N$. 
Furthermore 
$\widehat{M}\in\mathcal{S}(G[\widehat{N}])$ by Proposition~\ref{P1strong}. 
Given $v\in M$, we obtain $\{v\}\subsetneq\widehat{M}\subseteq\widehat{N}$. 
Since $\{v\}\in\Pi(G[\widehat{N}])$, $\widehat{M}=\widehat{N}$. 
Therefore $M=N$ because $M=\{v\in\widehat{M}:\{v\}\in\Pi(G[\widehat{M}])\}$ and 
$N=\{v\in\widehat{N}:\{v\}\in\Pi(G[\widehat{N}])\}$. 
\epr

Let $G$ be a graph such that $\max(\alpha_M(G),\omega_M(G))\geq 2$. 
Consider $M\in\mathbb{M}(G)$. 
By Proposition~\ref{P3th1}, $\lambda_G(\widehat{M})\in\{\medsquare,\filledmedsquare\}$ and 
$|\Pi_1(G[\widehat{M}])|=|M|\geq 2$. 
By Corollary~\ref{C2th1}, 
$$p(G)\leq q(G)\leq\max(\{\lceil\log_2(|\Pi_1(G[Y])|+1)\rceil:Y\in\mathcal{S}_{\geq 2}(G),\lambda_G(Y)\in\{\medsquare,\filledmedsquare\}\}).$$
We have also
\begin{gather}
\max(\{\lceil\log_2(|\Pi_1(G[Y])|+1)\rceil:Y\in\mathcal{S}_{\geq 2}(G),\lambda_G(Y)\in\{\medsquare,\filledmedsquare\}\})\notag\\
\shortparallel\tag{by Proposition~\ref{P3th1}}\\
\max(\{\lceil\log_2(|M|+1)\rceil:M\in\mathbb{M}(G)\})\notag\\
\shortparallel\notag\\
\lceil\log_2(\max(\alpha_M(G),\omega_M(G))+1)\rceil.\notag 
\end{gather}
Consequently 
\begin{gather}
p(G)\leq\lceil\log_2(\max(\alpha_M(G),\omega_M(G))+1)\rceil.\tag{Theorem~\ref{T1bound}}
\end{gather}

To obtain Corollary~\ref{C1bound}, we prove Lemma~\ref{L1bound}. 

~

\bprLonebound
Let $G$ be a graph such that $\max(\alpha_M(G),\omega_M(G))\geq 2$. 
There exists $S\in\mathcal{M}(G)$ such that 
$|S|=\max(\alpha_M(G),\omega_M(G))$ and $S$ is a clique or a stable set in $G$. 
Given an integer $p<\log_2(\max(\alpha_M(G),\omega_M(G)))$, 
consider any $p$-extension $H$ of $G$. 
We must prove that $H$ is not prime. 
We have $2^{\left|V(H)\setminus V(G)\right|}<\left|S\right|$ so that the function 
$S\longrightarrow 2^{V(H)\setminus V(G)}$, defined by $s\mapsto N_H(s)\cap (V(H)\setminus V(G))$, is not injective. 
There are $s\neq t\in S$ such that $v\longleftrightarrow_H\{s,t\}$ for every $v\in V(H)\setminus V(G)$. 
As $S$ is a module of $G$, we have $v\longleftrightarrow_H\{s,t\}$ 
for every $v\in V(G)\setminus S$. 
Since $S$ is a stable set in $G$, $\{s,t\}$ is a nontrival module of $H$. 
\epr

When a graph or its complement admits isolated vertices, we obtain the following. 

\begin{lem}\label{L2bound}
Given a graph $G$, if $\iota(G)\neq\emptyset$ or 
$\iota(\overline{G})\neq\emptyset$, then 
$$p(G)\geq\lceil\log_2(\max(\iota(G),\iota(\overline{G}))+1)\rceil.$$ 
\end{lem}

\bpr
By interchanging $G$ and $\overline{G}$, assume that $\iota(G)\geq\iota(\overline{G})$. 
Given $p<\lceil\log_2(\iota(G)+1)\rceil$, 
consider any $p$-extension $H$ of $G$.  
We have $2^{|V(H)\setminus V(G)|}\leq\iota(G)$ and we verify that $H$ is not prime. 

For each $v\in V(G)$ such that $N_G(v)=\emptyset$, we have $N_H(v)\subseteq V(H)\setminus V(G)$. 
Thus $(N_H)_{\restriction\{v\in V(G):N_G(v)=\emptyset\}}$ is a function from 
$\{v\in V(G):N_G(v)=\emptyset\}$ to 
$2^{V(H)\setminus V(G)}$. 
As observed in the proof of Lemma~\ref{L1stable}, if $(N_H)_{\restriction\{v\in V(G):N_G(v)=\emptyset\}}$ is not injective, then $\{u,v\}$ is a nontrivial module of $H$ when $u\neq v\in\{v\in V(G):N_G(v)=\emptyset\}$ with 
$N_H(u)=N_H(v)$. 
So assume that $(N_H)_{\restriction\{v\in V(G):N_G(v)=\emptyset\}}$ is injective. 
As $2^{|V(H)\setminus V(G)|}\leq\iota(G)$, we obtain that 
$(N_H)_{\restriction\{v\in V(G):N_G(v)=\emptyset\}}$ is bijective. 
Thus there is $u\in\{v\in V(G):N_G(v)=\emptyset\}$ such that $N_H(u)=\emptyset$, that is, 
$u\in\{v\in V(G):N_H(v)=\emptyset\}$. 
Therefore $H$ is not prime. 
\epr

The next is a simple consequence of Proposition~\ref{P3th1} which is useful in proving 
Theorem~\ref{T2bound}. 

\begin{cor}\label{C1th2}
Given a graph $G$ such that $\max(\alpha_M(G),\omega_M(G))\geq 2$, the elements of $\mathbb{M}(G)$ are pairwise disjoint. 
\end{cor}

\bpr
Consider $M,N\in\mathcal{M}_{{\rm max}}(G)$ such that $M\cap N\neq\emptyset$. 
Let $v\in M\cap N$. 
Since $\widehat{M},\widehat{N}\in\mathcal{S}(G)$ and $v\in\widehat{M}\cap\widehat{N}$, 
$\widehat{M}\subseteq\widehat{N}$ or $\widehat{N}\subseteq\widehat{M}$. 
For instance, assume that $\widehat{M}\subseteq\widehat{N}$. 
By Proposition~\ref{P1strong}, 
$\widehat{M}\in\mathcal{S}(G[\widehat{N}])$. 
Furthermore $\{v\}\in\Pi(G[\widehat{N}])$ by Proposition~\ref{P3th1}. 
As $\{v\}\subsetneq\widehat{M}\subseteq\widehat{N}$, we obtain 
$\widehat{M}=\widehat{N}$. 
Lastly, $M=\{w\in\widehat{M}:\{w\}\in\Pi(G[\widehat{M}])\}$ and 
$N=\{w\in\widehat{N}:\{w\}\in\Pi(G[\widehat{N}])\}$ by Proposition~\ref{P3th1}. 
Thus $M=N$. 
\epr

\section{Proof of Theorem~\ref{T2bound}}

Given a graph $G$, denote by $\mathbb{P}(G)$ the family of $M\in\mathcal{M}(G)$  such that 
$G[M]$ is prime. 
For every $M\in\mathbb{P}(G)$, $M\in\mathcal{S}(G)$ because $G[M]$ is prime. 
It follows that the elements of $\mathbb{P}(G)$ are pairwise disjoint. 
Thus the elements of $\mathbb{M}(G)\cup\mathbb{P}(G)$ are also by Corollary~\ref{C1th2}. 
Set $$I(G)=V(G)\setminus((\bigcup\mathbb{M}(G))\cup(\bigcup\mathbb{P}(G))).$$
We prove Theorem~\ref{T2bound} when $\max(\alpha_M(G),\omega_M(G))=2$. 

\begin{prop}\label{P1bound}
For every graph $G$ such that $\max(\alpha_M(G),\omega_M(G))=2$, 
\begin{equation*}
\text{$p(G)=2$ if and only if $\iota(G)=2$ or  $\iota(\overline{G})=2$.}
\end{equation*}
\end{prop}

\bpr 
It follows from Lemma~\ref{L1bound} and Theorem~\ref{T1bound} that $p(G)=1$ or 2. 
To begin, assume that $\iota(G)=2$ or  $\iota(\overline{G})=2$. 
By Lemma~\ref{L2bound}, $p(G)\geq 2$ and hence $p(G)=2$. 
Conversely, assume that $p(G)=2$. 
Let $a\not\in V(G)$. 
As $\max(\alpha_M(G),\omega_M(G))=2$, $|N|=2$ for each $N\in\mathbb{M}(G)$. 
Let $N_0\in\mathbb{M}(G)$. 
For $N\in\mathbb{P}(G)$, we have $G[N]$ is prime. 
By Lemma~\ref{L1extension}, $G[N]$ admits a prime extension $H_N$ defined on $N\cup\{a\}$. 
We consider any 1-extension $H$ of $G$ to $V(G)\cup\{a\}$ satisfying the following. 
\begin{enumerate}
\item For each $N\in\mathbb{M}(G)$, $a\not\longleftrightarrow_H N$.
\item For each $N\in\mathbb{P}(G)$, $H[N\cup\{a\}]=H_N$. 
\item Let $v\in I(G)$. There is $i\in\{0,1\}$ such that $(v,N_0)_G=i$. 
We require that $(v,a)_H\neq i$. 
\end{enumerate}
To begin, we prove that $\mathcal{S}_{\geq 2}(G)\cap\mathcal{M}(H)=\emptyset$. 
Given $M\in\mathcal{S}_{\geq 2}(G)$, we have to verify that 
$a\not\longleftrightarrow_H M$. 
Let $N$ be a minimal element under inclusion of $\{N'\in\mathcal{S}_{\geq 2}(G):N'\subseteq M\}$. 
By Proposition~\ref{P1strong}, $\Pi(G[N])\subseteq\mathcal{S}(G)$. 
By minimality of $N$, $\Pi(G[N])=\Pi_1(G[N])$ so that $G[N]$ and  $G[N]/\Pi(G[N])$ are isomorphic by Proposition~\ref{P1quotient}.1. 
We distinguish the following two cases. 
\begin{itemize}
\item Assume that $\lambda_G(N)=\sqcup$. 
We obtain that $G[N]$ is prime, that is, $N\in\mathbb{P}(G)$. 
As $H[N\cup\{a\}]$ is prime, $a\not\longleftrightarrow_H N$. 
\item Assume that $\lambda_G(N)\in\{\medsquare,\filledmedsquare\}$. 
By Proposition~\ref{P3th1}, $N\in\mathbb{M}(G)$. 
Thus $|N|=2$ and $a\not\longleftrightarrow_H N$ by definition of $H$.
\end{itemize}
In both cases, $a\not\longleftrightarrow_H N$ and hence 
$a\not\longleftrightarrow_H M$. 

Now we prove that $\mathcal{M}_{\geq 2}(G)\cap\mathcal{M}(H)=\emptyset$. 
Let $M\in\mathcal{M}_{\geq 2}(G)$. 
Since $\mathcal{S}_{\geq 2}(G)\cap\mathcal{M}(H)=\emptyset$, assume that 
$M\not\in\mathcal{S}_{\geq 2}(G)$. 
Set $Q=\{X\in\Pi(G[\widehat{M}]):X\cap M\neq\emptyset\}$. 
By Proposition~\ref{P1module}.1, $M\in\mathcal{M}(G[\widehat{M}])$. 
By definition of $\widehat{M}$, $|Q|\geq 2$. 
Thus $M=\bigcup Q$ because $\Pi(G[\widehat{M}])\subseteq\mathcal{S}(G[\widehat{M}])$. 
Furthermore $Q\neq\Pi(G[\widehat{M}])$ because $M\not\in\mathcal{S}_{\geq 2}(G)$. 
By Proposition~\ref{P1quotient}.2, 
$Q\in\mathcal{M}(G[\widehat{M}]/\Pi(G[\widehat{M}]))$. 
As $2\leq|Q|<|\Pi(G[\widehat{M}])|$, 
$\lambda_G(\widehat{M})\in\{\medsquare,\filledmedsquare\}$. 
If there is $X\in Q\cap\Pi_{\geq 2}(G[\widehat{M}])$, then 
$a\not\longleftrightarrow_H X$ by what precedes and hence 
$a\not\longleftrightarrow_H M$. 
Assume that $Q\subseteq\Pi_1(G[\widehat{M}])$. 
We obtain that $M$ is a clique or a stable set in $G$. 
Since $\max(\alpha_M(G),\omega_M(G))=2$, $M\in\mathbb{M}(G)$ and 
$a\not\longleftrightarrow_H M$ by definition of $H$. 

As $p(G)=2$, $H$ admits a nontrivial module $M_H$. 
We have $a\in M_H$ because $\mathcal{M}_{\geq 2}(G)\cap\mathcal{M}(H)=\emptyset$. 

First, we show that $N\subseteq M_H$ for each $N\in\mathbb{P}(G)$. 
By Proposition~\ref{P1module}.1, 
$M_H\cap(N\cup\{a\})\in\mathcal{M}(H[N\cup\{a\}])$. 
Since $H[N\cup\{a\}]$ is prime and $a\in M_H\cap (N\cup\{a\})$, we obtain either 
$(M_H\setminus\{a\})\cap N=\emptyset$ or 
$N\subseteq M_H\setminus\{a\}$. 
Suppose for a contradiction that $(M_H\setminus\{a\})\cap N=\emptyset$. 
By Proposition~\ref{P1module}.1, $M_H\setminus\{a\}\in\mathcal{M}(G)$. 
There is $i\in\{0,1\}$ such that $(M_H\setminus\{a\},N)_G=i$ by Proposition~\ref{P1module}.3. 
Therefore $(a,N)_H=i$ which contradicts the fact that 
$H[N\cup\{a\}]$ is prime. 
It follows that $N\subseteq M_H$. 
Thus 
\begin{equation}\label{E1P1bound}
\bigcup\mathbb{P}(G)\subseteq M_H.
\end{equation}

Second, we show that $N\cap M_H\neq\emptyset$ for each 
$N\in\mathbb{M}(G)$. 
Otherwise consider $N\in\mathbb{M}(G)$ such that $N\cap M_H=\emptyset$. 
There is $i\in\{0,1\}$ such that $(M_H\setminus\{a\},N)_G=i$. 
Thus $(a,N)_H=i$ which contradicts $a\not\longleftrightarrow_H N$. 
Therefore 
\begin{equation}\label{E2P1bound}
N\cap M_H\neq\emptyset\quad\text{for each}\quad N\in\mathbb{M}(G).
\end{equation}

Third, let  $v\in I(G)$. 
By \eqref{E2P1bound}, $N_0\cap M_H\neq\emptyset$. 
Since $(v,N_0\cap M_H)_G\neq (v,a)_H$, $v\in M_H$. 
Hence 
\begin{equation}\label{E3P1bound}
I(G)\subseteq M_H.
\end{equation}
By \eqref{E1P1bound} and \eqref{E3P1bound}, 
\begin{equation}\label{E4P1bound}
V(G)\setminus M_H\subseteq\mathbb{M}(G).
\end{equation}

To conclude, consider $v\in V(H)\setminus M_H$. 
By \eqref{E4P1bound}, there is $N_v\in\mathbb{M}(G)$ such that $v\in N_v$. 
By interchanging $G$ and $\overline{G}$, assume that $N_v$ is a stable set in $G$. 
Since $v\longleftrightarrow_HM_H $ and $(v,N_v\cap M_H)_G=0$, we obtain $(v,M_H)_H=0$. 
Let $N\in\mathbb{M}(G)\setminus\{N_v\}$. 
By Corollary~\ref{C1th2}, $N\cap N_v=\emptyset$. 
As $N\cap M_H\neq\emptyset$ by \eqref{E2P1bound}, we have $(v,N\cap M_H)_G=0$ and hence $(v,N)_G=0$. 
It follows that $N_G(v)=\emptyset$. 
Therefore $(N_v,V(G)\setminus N_v)_G=0$ because $N_v\in\mathcal{M}(G)$. 
Since $N_v$ is a stable set in $G$, we obtain $N_v\subseteq\{u\in V(G):N_G(u)=\emptyset\}$. 
Clearly $\{u\in V(G):N_G(u)=\emptyset\}\in\mathcal{M}(G)$ and $\{u\in V(G):N_G(u)=\emptyset\}$ is a stable set in $G$. 
Thus $\iota(G)\leq\max(\alpha_M(G),\omega_M(G))=2$. 
Consequently $N_v=\{u\in V(G):N_G(u)=\emptyset\}$. 
\epr

\bprthtwobound
Consider a graph $G$ such that $\max(\alpha_M(G),\omega_M(G))=2^k$ where $k\geq 1$. 
It follows from Lemma~\ref{L1bound} and Theorem~\ref{T1bound} that $p(G)=k$ or $k+1$. 
To begin, assume that $\iota(G)=2^k$ or $\iota(\overline{G})=2^k$. 
By Lemma~\ref{L2bound}, $p(G)\geq k+1$ and hence $p(G)=k+1$. 

Conversely, assume that $p(G)=k+1$. 
If $k=1$, then it suffices to apply Proposition~\ref{P1bound}. 
Assume that $k\geq 2$. 
For convenience set 
$$\mathbb{M}_{{\rm max}}(G)=\{N\in\mathbb{M}(G):|N|=\max(\alpha_M(G),\omega_M(G))\}.$$
With each $N\in\mathbb{M}_{{\rm max}}(G)$ associate 
$w_N\in N$. 
Set $W=\{w_N:N\in\mathbb{M}_{{\rm max}}(G)\}$. 

We prove that $\max(\alpha_M(G-W),\omega_M(G-W))=2^k-1$. 
Let $N\in\mathbb{M}_{{\rm max}}(G)$. 
By Corollary~\ref{C1th2}, the elements of $\mathbb{M}_{{\rm max}}(G)$ are pairwise disjoint. 
Thus $N\setminus W=N\setminus\{w_N\}$. 
Clearly $N\setminus\{w_N\}$ is a clique or a stable set in $G-W$. 
Furthermore $N\setminus\{w_N\}\in\mathcal{M}(G-W)$. 
Therefore $2^k-1=|N\setminus\{w_N\}|\leq\max(\alpha_M(G-W),\omega_M(G-W))$. 
Now consider $N'\in\mathbb{M}_{{\rm max}}(G-W)$. 
We show that $N'\in\mathcal{M}(G)$. 
We have to verify that for each $N\in\mathbb{M}_{{\rm max}}(G)$, 
$w_N\longleftrightarrow_G N'$. 
Let $N\in\mathbb{M}_{{\rm max}}(G)$. 
First, asume that there is $v\in(N\setminus\{w_N\})\setminus N'$. 
We have $v\longleftrightarrow_G N'$. 
As $N$ is a clique or a stable set in $G$, $\{v,w_N\}\in\mathcal{M}(G[N])$. 
By Proposition~\ref{P1module}.2, $\{v,w_N\}\in\mathcal{M}(G)$. 
Thus $w_N\longleftrightarrow_G N'$. 
Second, assume that $N\setminus\{w_N\}\subseteq N'$. 
Clearly $w_N\longleftrightarrow_G N'$ when $N\setminus\{w_N\}=N'$. 
Assume that $N'\setminus (N\setminus\{w_N\})\neq\emptyset$. 
By interchanging $G$ and $\overline{G}$, assume that $N'$ is a clique in $G-W$. 
As $N\setminus\{w_N\}\subseteq N'$ and $|N\setminus\{w_N\}|\geq 2$, we obtain that $N$ is a clique in $G$. 
Since $(N\setminus\{w_N\},N'\setminus N)_G=1$ and since $N\in\mathcal{M}(G)$, we have 
$(w_N,N'\setminus N)_G=1$. 
Furthermore $(w_N,N\setminus\{w_N\})_G=1$ because $N$ is a clique in $G$. 
Therefore $(w_N,N')_G=1$. 
Consequently $N'\in\mathcal{M}(G)$. 
As $N'$ is a clique in $G$, there is $M\in\mathbb{M}(G)$ such that $M\supseteq N'$. 
If $M\not\in\mathbb{M}_{{\rm max}}(G)$, then 
$|N'|\leq|M|<\max(\alpha_M(G),\omega_M(G))$. 
If $M\in\mathbb{M}_{{\rm max}}(G)$, then 
$N'\subseteq M\setminus\{w_M\}$ and hence $|N'|<|M|=\max(\alpha_M(G),\omega_M(G))$. 
In both cases, we have $|N'|=\max(\alpha_M(G-W),\omega_M(G-W))<\max(\alpha_M(G),\omega_M(G))$. 
It follows that $\max(\alpha_M(G-W),\omega_M(G-W))=2^k-1$. 

By Lemma~\ref{L1bound} and Theorem~\ref{T1bound}, $p(G-W)=k$ and hence there exists a prime $k$-extension $H'$ of $G-W$. 
We extend $H'$ to $V(H')\cup W$ as follows. 
Let $N\in\mathbb{M}_{{\rm max}}(G)$. 
Consider the function $f_N:N\setminus\{w_N\}\longrightarrow 2^{V(H')\setminus V(G-W)}$ defined by $v\mapsto N_{H'}(v)\setminus V(G-W)$ for $v\in N\setminus\{w_N\}$. 
Since $H'$ is prime, $f_N$ is injective. 
As $|N\setminus\{w_N\}|=2^k-1$ and $|2^{V(H')\setminus V(G-W)}|=2^k$, there is a unique 
$X_N\subseteq V(H')\setminus V(G-W)$ such that $f_N(v)\neq X_N$ for every 
$v\in N\setminus\{w_N\}$. 
Let $H$ be the extension of $H'$ to $V(H')\cup W$ such that 
$N_H(w_N)\cap (V(H')\setminus V(G-W))=X_N$ for each 
$N\in\mathbb{M}_{{\rm max}}(G)$. 
As $p(G)=k+1$, $H$ is not prime. 
Consider a nontrivial module $M_H$ of $H$. 

Observe the following. 
Given $N\neq N'\in\mathbb{M}_{{\rm max}}(G)$,
\begin{equation}\label{E1th2}
\left.\begin{array}{r}
N\cap M_H\neq\emptyset\\
\text{and}\\
N'\cap M_H\neq\emptyset\\
\end{array}\right\}\ \Longrightarrow M_H\supseteq V(H').
\end{equation}
Indeed, by Proposition~\ref{P1module}.1, $M_H\cap V(G)\in\mathcal{M}(G)$. 
Since $\widehat{N},\widehat{N'}\in\mathcal{S}(G)$ and since 
$(M_H\cap V(G))\cap\widehat{N}\neq\emptyset$ and 
$(M_H\cap V(G))\cap\widehat{N'}\neq\emptyset$, $M_H\cap V(G)$ is comparable to 
$\widehat{N}$ and $\widehat{N'}$ under inclusion. 
Suppose for a contradiction that $M_H\cap V(G)\subsetneq\widehat{N}$ and 
$M_H\cap V(G)\subsetneq\widehat{N'}$. 
It follows that $N'\cap\widehat{N}\neq\emptyset$ and 
$N\cap\widehat{N'}\neq\emptyset$. 
As $\widehat{N'}\in\mathcal{S}(G)$, $\widehat{N'}\subsetneq N$ or $N\subseteq\widehat{N'}$. 
In the first instance, it follows from Proposition~\ref{P1strong} that 
$\widehat{N'}$ would be a nontrivial strong module of $G[N]$ which contradicts the fact that $N$ is a clique or a stable set in $G$. 
Thus $N\subseteq\widehat{N'}$ and hence $\widehat{N}\subseteq\widehat{N'}$. 
Similarly $N'\subseteq\widehat{N}$ and $\widehat{N'}\subseteq\widehat{N}$. 
Therefore $\widehat{N}=\widehat{N'}$ and it would follow from Proposition~\ref{P3th1} that $N=N'$. 
Consequently $\widehat{N}\subseteq(M_H\cap V(G))$ or 
$\widehat{N'}\subseteq(M_H\cap V(G))$. 
For instance, assume that $\widehat{N}\subseteq(M_H\cap V(G))$. 
By Proposition~\ref{P1module}.1, $M_H\cap V(H')\in\mathcal{M}(H')$. 
Furthermore $(M_H\cap V(H'))\supseteq (N\setminus W)$ and 
$N\setminus W=N\setminus\{w_N\}$ by Corollary~\ref{C1th2}. 
Since $H'$ is prime, we have $V(H')\subseteq M_H$. 
It follows that \eqref{E1th2} holds. 

As $H'$ is prime and $M_H\cap V(H')\in\mathcal{M}(H')$, we have 
either $|M_H\cap V(H')|\leq 1$ or $M_H\supseteq V(H')$. 
For a contradiction, suppose that $|M_H\cap V(H')|\leq 1$. 
There is $N\in\mathbb{M}_{{\rm max}}(G)$ such that 
$w_N\in M_H$. 
It follows from \eqref{E1th2} that 
\begin{equation}\label{E2th2}
N'\cap M_H=\emptyset\ \text{for each 
$N'\in\mathbb{M}_{{\rm max}}(G)\setminus\{N\}$.}
\end{equation} 
Thus $M_H\cap W=\{w_N\}$ and there is $v\in V(H')$ such that $M_H\cap V(H')=\{v\}$. 
Clearly $M_H=\{v,w_N\}$ and we distinguish the following two cases to obtain a contradiction.
\begin{itemize}
\item Suppose that $v\in V(G-W)$. 
By Proposition~\ref{P1module}.1, $\{v,w_N\}\in\mathcal{M}(G)$. 
Therefore there is $N'\in\mathbb{M}_{{\rm max}}(G)$ such that $N'\supseteq\{v,w_N\}$. 
By \eqref{E2th2}, $N=N'$ and we would obtain $N_H(w_N)\cap (V(H')\setminus V(G-W))=f_N(v)$. 
\item Suppose that $v\in V(H')\setminus V(G-W)$. 
There is $i\in\{0,1\}$ such that $(w_N,N\setminus\{w_N\})_G=i$. 
We obtain $(v,N\setminus\{w_N\})_{H'}=i$ because $\{v,w_N\}\in\mathcal{M}(H)$. 
Since $f_N$ is injective, the function $g_N:N\setminus\{w_N\}\longrightarrow 2^{((V(H')\setminus V(G-W))\setminus\{v\})}$, defined by 
$g_N(u)=f_N(u)\setminus\{v\}$ for $u\in N\setminus\{w_N\}$, is injective as well. 
We would obtain $2^k-1\leq 2^{k-1}$. 
\end{itemize}
Consequently $V(H')\subseteq M_H$. 
As $M_H$ is a nontrivial module of $H$, there exists $N\in\mathbb{M}_{{\rm max}}(G)$ such that $w_N\not\in M$. 
By interchanging $G$ and $\overline{G}$, assume that $N$ is a stable set in $G$. 
We have $(w_N,N\setminus\{w_N\})_G=0$ and hence $(w_N,V(H'))_H=0$. 
In particular $(w_N,V(G-W))_G=0$. 
Given $N'\in\mathbb{M}_{{\rm max}}(G)\setminus\{N\}$, we obtain 
$(w_N,N'\setminus\{w_{N'}\})_G=0$. 
Since $N'\in\mathcal{M}(G)$, $(w_N,w_{N'})_G=0$. 
It follows that $N_G(w_N)=\emptyset$. 
As at the end of the proof of Proposition~\ref{P1bound}, we conclude by 
$N=\{u\in V(G):N_G(u)=\emptyset\}$. 
\epr

Lastly, we examine the graphs $G$ such that $\alpha_M(G)=\omega_M(G)=1$. For these, 
$\mathbb{M}(G)=\emptyset$. Thus either $|V(G)|\leq 1$ or $|V(G)|\geq 4$ and $G$ is not prime. 

\begin{prop}\label{m(G)=1}
For every non prime graph $G$ such that $|V(G)|\geq 4$ and $\alpha_M(G)=\omega_M(G)=1$, we have $p(G)=1$. 
\end{prop}

\bpr
Consider a minimal element $N_{{\rm min}}$ of $\mathcal{S}_{\geq 2}(G)$. 
By Proposition~\ref{P1strong}, 
$\Pi($ $G[N_{{\rm min}}])\subseteq\mathcal{S}(G)$. 
By minimality of $N_{{\rm min}}$, $\Pi(G[N_{{\rm min}}])=\Pi_1(G[N_{{\rm min}}])$. 
Thus $G[N_{{\rm min}}]$ and $G[N_{{\rm min}}]/\Pi(G[N_{{\rm min}}])$ are isomorphic by Proposition~\ref{P1quotient}.1. 
If $\lambda_G($ $N_{{\rm min}})\in\{\medsquare,\filledmedsquare\}$, then $N_{{\rm min}}$ is a clique or a stable set in $G$ and there would be $N\in\mathbb{M}(G)$ such that 
$N\supseteq N_{{\rm min}}$. 
Therefore $\lambda_G(N_{{\rm min}})=\sqcup$ and $N_{{\rm min}}\in\mathbb{P}(G)$. 

Let $a\not\in V(G)$. 
For each $N\in\mathbb{P}(G)$, $G[N]$ is prime. 
By Lemma~\ref{L1extension}, $G[N]$ admits a prime 1-extension $H_N$ to 
 $N\cup\{a\}$. 
We consider the 1-extension $H$ of $G$ to $V(G)\cup\{a\}$ satisfying the following. 
\begin{enumerate}
\item For each $N\in\mathbb{P}(G)$, $H[N\cup\{a\}]=H_N$. 
\item Let $v\in I(G)$. 
There is $i\in\{0,1\}$ such that $(v,N_{{\rm min}})_G=i$. 
We require that $(v,a)_H\neq i$. 
\end{enumerate}

We proceed as in the proof of Proposition~\ref{P1bound}, to show that 
$\mathcal{M}_{\geq 2}(G)\cap\mathcal{M}(H)=\emptyset$. 
To begin, we prove that $\mathcal{S}_{\geq 2}(G)\cap\mathcal{M}(H)=\emptyset$. 
Given $M\in\mathcal{S}_{\geq 2}(G)$, we have to verify that 
$a\not\longleftrightarrow_H M$. 
Let $N$ be a minimal element under inclusion of $\{N'\in\mathcal{S}_{\geq 2}(G):N'\subseteq M\}$. 
We obtain that $\Pi(G[N])=\Pi_1(G[N])$ so that $G[N]$ and  $G[N]/\Pi(G[N])$ are isomorphic by Proposition~\ref{P1quotient}.1. 
If $\lambda_G(N)\in\{\medsquare,\filledmedsquare\}$, then $N$ is a clique or a stable set in $G$ and there would be $N'\in\mathbb{M}(G)$ such that 
$N'\supseteq N$. 
Thus $\lambda_G(N)=\sqcup$. 
We obtain that $G[N]$ is prime, that is, $N\in\mathbb{P}(G)$. 
Since $H[N\cup\{a\}]$ is prime, $a\not\longleftrightarrow_H N$ and hence 
$a\not\longleftrightarrow_H M$. 

Now we prove that $\mathcal{M}_{\geq 2}(G)\cap\mathcal{M}(H)=\emptyset$. 
Let $M\in\mathcal{M}_{\geq 2}(G)$. 
Since $\mathcal{S}_{\geq 2}(G)\cap\mathcal{M}(H)=\emptyset$, assume that 
$M\not\in\mathcal{S}_{\geq 2}(G)$. 
Set $Q=\{X\in\Pi(G[\widehat{M}]):X\cap M\neq\emptyset\}$. 
We obtain that $M=\bigcup Q$, $|Q|\geq 2$ and 
$\lambda_G(\widehat{M})\in\{\medsquare,\filledmedsquare\}$. 
If $|\Pi_1(G[\widehat{M}])|\geq 2$, then we would have 
$\{v\in\widehat{M}:\{v\}\in\Pi(G[\widehat{M}])\}\in\mathbb{M}(G)$ by Proposition~\ref{P3th1}. 
Consequently $|\Pi_1(G[\widehat{M}])|\leq 1$ and there is 
$X\in Q\cap\Pi_{\geq 2}(G[\widehat{M}])$. 
By what precedes $a\not\longleftrightarrow_H X$ and hence 
$a\not\longleftrightarrow_H M$. 

Lastly, we establish that $H$ is prime. 
Let $M_H\in\mathcal{M}_{\geq 2}(H)$. 
As previously shown, $a\in M$. 
We show that $N\subseteq M_H$ for each $N\in\mathbb{P}(G)$. 
By Proposition~\ref{P1module}.1, 
$M_H\cap(N\cup\{a\})\in\mathcal{M}(H[N\cup\{a\}])$. 
Since $H[N\cup\{a\}]$ is prime and $a\in M_H\cap (N\cup\{a\})$, we obtain either 
$(M_H\setminus\{a\})\cap N=\emptyset$ or 
$N\subseteq M_H\setminus\{a\}$. 
Suppose for a contradiction that $(M_H\setminus\{a\})\cap N=\emptyset$. 
By Proposition~\ref{P1module}.1, $M_H\setminus\{a\}\in\mathcal{M}(G)$. 
There is $i\in\{0,1\}$ such that $(M_H\setminus\{a\},N)_G=i$ by Proposition~\ref{P1module}.3. 
Therefore $(a,N)_H=i$ which contradicts the fact that 
$H[N\cup\{a\}]$ is prime. 
It follows that $N\subseteq M_H$ for each $N\in\mathbb{P}(G)$. 
In particular $N_{{\rm min}}\subseteq M$. 
Let $v\in I(G)$. 
As $(v,N_{{\rm min}})_G\neq (v,a)_H$, $v\in M_H$. 
Consequently $M_H=V(H)$. 
\epr

\end{document}